\documentclass[12pt]{amsart}
\usepackage{amssymb,amsmath,amsthm,amsfonts,mathrsfs,mathtools}

\usepackage{a4wide}
\usepackage[utf8]{inputenc}
\usepackage[T1]{fontenc}
\usepackage{color}
\usepackage{esint}

\newcommand{\R}{\mathbb{R}}

\newcommand{\intbar}{\fint}

\newcommand{\N}{\mathbb{N}}

\DeclareMathOperator{\dist}{dist}

\usepackage{graphicx}
\newcommand{\mres}{\mathbin{\vrule height 1.6ex depth 0pt width
		0.13ex\vrule height 0.13ex depth 0pt width 1.3ex}} 

\newtheorem{theorem}{Theorem}[section]
\newtheorem{proposition}{Proposition}[section]

\newtheorem{lemma}{Lemma}[section]
\newtheorem{corollary}{Corollary}[section]
\newtheorem{definition}{Definition}[section]
\newtheorem{remark}{Remark}[section]

\numberwithin{equation}{section} \numberwithin{theorem}{section}
\numberwithin{proposition}{section} \numberwithin{lemma}{section}
\numberwithin{corollary}{section}
\numberwithin{definition}{section} \numberwithin{remark}{section}

\let\O=\Omega

\let\eps=\varepsilon

\begin{document}

\title{A field-road system  with  a rectifiable set}

\author{Matthieu Bonnivard, Romain Ducasse, Antoine Lemenant, Alessandro Zilio}

\begin{abstract} The aim of this paper is to define a  field-road system in 2D where   the road is a  merely 1D-rectifiable set. For this purpose we introduce a general setting in order to define a parabolic problem onto a rectifiable set, which is coupled with another more classical parabolic problem outside this set, with transmission conditions.    
\end{abstract}

\maketitle

\setcounter{tocdepth}{1}

\tableofcontents

\section{Introduction}

The so-called \emph{road-field} system introduced by Berestycki,  Roquejoffre, and  Rossi in \cite{BRR} aims to study the impact of a road in the propagation of a species or a disease. In this model, the species usually thrives in a domain $\Omega\subset \R^2$ with moderate diffusion capacities. This domain is crossed by a road $K\subset \Omega$, where a faster diffusion occurs. The questions addressed by Berestycki \emph{et al}.\! include the identification of the road's effect on the overall invasion speed of the species, as well as its precise quantification.

In~\cite{BRR}, the road $K$  is  assumed to be a straight segment and in the linear setting the system contains  reaction-diffusion type equations that are coupled with a transmission condition as follows:
\begin{equation}\label{systeme1}
\left\{
\begin{array}{rll}
\partial_t v - a \Delta v &= 0   & \text{ in }\Omega\backslash K,\\
\partial_t u - b \Delta_K u &=  \nu v\vert_K - \mu u& \text{ on } K,\\
a( \partial_n^+ v + \partial_n^- v)\vert_K&  =  \mu u - \nu v\vert_K  & \text{ on }   K,
\end{array}
\right.
\end{equation}
with additional standard boundary conditions. \textcolor{black}{Here $( \partial_n^+ v + \partial_n^- v)\vert_K$ stands for the sum of the outward derivatives of $v$ on $K$ coming from the two sides of the segment $K$.}

A possible interpretation for the system is the following. We have a population, whose density on $\O$ is given by $v$ and whose density on $K$ is given by $u$. The individuals diffuse on $\O$ and $K$, with diffusivity constant $a, b>0$ respectively. 

The individuals located in the field by the road (whose density is $v\vert_K$, that is the trace of $v$ on $K$) enter the road with a rate $\nu$.

The individuals on the road leave it (and enter the field) with a rate $\mu$. The road has two sides, we assume that the individuals have the same probability to leave to one side or the other. In the system, we denote $\partial_n^+, \partial_n^-$ the "upward" and "downward" flux. When the road is a straight line, this is easily defined.

When the exchanges coefficients satisfy $\mu = \nu =0$ for instance, the problem becomes ``uncoupled", that is, individuals on the field do not see the road anymore and diffuse as if it were not here (this case is somewhat degenerate).

\medskip

  The aim of this paper is to generalize the previous model to the case where the road $K$ is allowed to be a general network composed by curves, of even more generally, a compact connected set $K$ with $\mathcal{H}^1(K)<+\infty$. In such setting, since the network $K$ can be very wild, it is not clear whether or not we can understand this system in a classical sense. However, we can define solutions in  a weak sense by use of a certain Sobolev space $H^1(K)$ associated to the rectifiable set~$K$.

  One of the main contributions of this paper is indeed a new approach to define a proper Sobolev space $H^1(K)$ associated to a 1D-rectifiable set $K$, and its associated Dirichlet energy. This space enjoys certain natural properties such as the compact embedding into the space $L^2(K,d\mathcal{H}^1)$. Actually, our space $H^1(K)$ coincides with the space $H^1_\mu$ defined in \cite{BBS}, with $\mu=\mathcal{H}^1|_K$. The approach in \cite{BBS} is very general and uses powerful abstract tools such as convex duality. In contrast, thanks to the fact that $K$ is connected with finite $\mathcal{H}^1$ measure, we propose a much simpler approach which appears more efficient to define the weak formulation for a parabolic problem.  The construction of this space is done in Section \ref{sectionH^1}. {\color{black} Let us mention that an alternative approach can be found in \cite{ABSO}, Section 2.9. The difference is that in \cite{ABSO} the authors use a more embedded point of view invoking a special parameterization given by  \cite{alberti}. Our approach is different, more extrinsic, in the sense that it does not rely on a given parameterization.  However, since in  \cite{ABSO} the functions are  moreover assumed to be globally continuous, we expect both constructions to produce the same result.}

  Once the space $H^1(K)$ is well defined, we can define a weak solution for the parabolic system \eqref{systeme1} and prove the existence of a solution. This is done in Section ~\ref{Existence}. The method relies on a suitable spectral decomposition performed in Section~\ref{Spectral}, which also allows us to derive the long time behaviour of the weak solution (see Corollary~\ref{Cor:longtime}). We conclude the paper by comparing in Section~\ref{Sect:conclusion} this asymptotic behaviour to the one observed in the absence of a road, which enables us to identify a quantitative criterion measuring the effect of a road on the spreading of a population in this context.

\section{Definition of the Sobolev space $H^1(K)$ and basic properties}

\label{sectionH^1} 
We denote by $\mathcal{H}^1$ the one dimensional Hausdorff measure in $\R^2$. In what follows we denote by $L^2(K)$ the usual complete space  $L^2(K,d\mathcal{H}^1)$ containing measurable functions $u$ such that $\int_{K} u^2 \, d\mathcal{H}^1<+\infty$, and defined $\mathcal{H}^1$-a.e.\! on $K$. In what follows we describe a general strategy in order  to define a Dirichlet energy associated to the 1-rectifiable set $K$.  
 
 \subsection{Definition}
 \label{defH^1section}
 Let $\Omega \subset \R^2$ be a smooth open set. We will denote by $\mathcal{K}(\Omega)$ the class of all  $K\subset \overline\Omega$ being compact, connected, and satisfying $\mathcal{H}^1(K)<+\infty$. The following well-known facts are standard.
 
 \begin{proposition}\cite[Proposition 30.1 p.186]{d} \label{lenghtspace} 
 Let $K\in \mathcal{K}(\Omega)$. Then
\begin{itemize}
\item $K$ is   a $1$-rectifiable set.
\item $K$ is arcwise connected: for every $x,y \in K$ there exists an injective Lipschitz curve $\gamma:[0,1]\to K$ such that $\gamma(0)=x$, $\gamma(1)=y$, and $\mathcal{H}^1(\gamma([0,1]))=\dist_K(x,y)$, where $\dist_K(x,y)$ is the geodesic distance in $K$, defined by
\begin{eqnarray}
\dist_K(x,y):=\inf \left\{  \int_{0}^1 |\gamma'(t)|dt \; \big|\;  \gamma\in Lip([0,1],K ),   \gamma(0)= x,  \gamma(1)=y \right\}. \label{defdK}
\end{eqnarray}
\end{itemize} 
  \end{proposition}

 In particular, for   $K\in \mathcal{K}(\Omega)$ we know the existence of  an approximative tangent line at $\mathcal{H}^1$-a.e.\! point $x\in K$, and  for any such point $x$ we can choose a unit tangent vector $\tau_K(x)$ in direction of that line. \textcolor{black}{Observe that, if $K$ has a locally smooth parametrization $\gamma\in C^1([0,1])$, then $\tau_K(x)$ is a unit vector parallel to the tangent vector $\gamma'$.}

For every smooth function $u \in C^\infty(\R^2)$ we introduce 
$$N(u)=\left(\int_{K} u^2(y) \; d\mathcal{H}^1(y) + \int_{K} |\nabla u \cdot \tau_K(y)|^2 \; d\mathcal{H}^1(y)\right)^\frac{1}{2}.$$
Notice that $N(u)$ involves only the trace on $K$ of the smooth function $u$ defined on the whole $\R^2$. Then we consider the space $\mathcal{D}(K)$ as being the restriction on $K$ of  $C^\infty(\R^2)$ functions, and we endow this space with the norm  $N$. In particular,  $\mathcal{D}(K)$  is a subspace of $L^2(K)$.  Finally, we define $H^1(K)$ as follows. 

\begin{definition}[Space $H^1(K)$ and $\nabla_K u$] For $K\in\mathcal{K}(\Omega)$, we define $H^1(K)$ as the completion of $\mathcal{D}(K)$ for the norm $N$. In particular, $H^1(K)$ is a closed subspace of $L^2(K)$ for which $\mathcal{D}(K)$ is a dense subset. 

We define similarly the space $H^1_{0,\partial \Omega}(K)$ as the subspace of $H^1(K)$ of functions that vanish on $\partial \Omega \cap  K$, that is, the completion of $\{u \in \mathcal{D}(K) \ : \ u = 0 \ \text{on}\ \partial \Omega \cap  K \}$. If $K\cap \partial\Omega = \emptyset$, then $H^1_{0,\partial \Omega}(K)=H^1(K)$.

For any $u \in H^1(K)$ we define $\nabla_K u \in L^2(K;\R^2)$ as the $L^2$~limit of the projection (a.e.) of $\nabla u_n$ on $K$, that is, $(\nabla u_n\cdot \tau_K)\tau_K$, for $u_n \to u$ in the norm $N$. In particular $\nabla_K u$ does not depend on the choice of the sequence $u_n \in \mathcal{D}(K)$  such that $u_n \to u$ in $H^1(K)$, in the equivalent class of Cauchy sequences for the norm $N$, which justifies the definition.

\end{definition}

The construction of $H^1(K)$ is rather standard. A way to define it rigorously  is for instance by considering all Cauchy sequences for $N$ in $\mathcal{D}(K)$, on which one defines the following equivalence relation: two Cauchy sequences $u_n$ and $v_n$ are equivalent if and only if $N(u_n-v_n)\to 0$. Then $H^1(K)$ is the quotient of all Cauchy sequences by this relation. It is easy to see  that this space is a complete space for which $\mathcal{D}(K)$ is a dense subset.

  In particular, we can consider  $u \in H^1(K)$ as being a function $u\in L^2(K)$ for which there exists a sequence $u_n \in C^\infty(\R^2)$ such that $u_n\vert_K \to u$ in $L^2$ and   $(\nabla u_n \cdot \tau_K)\tau_K$ has a limit in $L^2(K;\R^2)$. For $u \in H^1(K)$ we will denote by $\nabla_K u$ the  limit  of $(\nabla u_n \cdot \tau_K)\tau_K$. By construction,  the limit of $(\nabla u_n \cdot \tau_K)\tau_K$ does not depend on the choice of the sequence $u_n$, chosen in the equivalent class of Cauchy sequences, and $\nabla_K u$ is therefore well defined.

Next, our Dirichlet energy, defined on $H^1(K)$, is given by
$$\int_{K} |\nabla_K u|^2 \;d\mathcal{H}^1.$$

Of course if $u\in  C^\infty(\R^2)$, $\nabla_K u = (\nabla u \cdot \tau_K)\tau_K$ $\mathcal{H}^1$-a.e., thus the Dirichlet energy coincides with the natural one in that case.  

It is also easy to see that our definition coincides with the one of \cite{BBS} for the particular case $\mu=\mathcal{H}^1|_K$.

\subsection{Compact embedding of $H^1(K)$ in $L^2(K)$.}

Now we would like to establish an analogue of Sobolev embedding and Rellich theorem within our context.  

\begin{proposition}\label{embb} If $K \in \mathcal{K}(\Omega)$, then for all $u\in H^1(K)$,  
\begin{eqnarray}
|u(x)-u(y)|\leq \dist_K(x,y)^\frac{1}{2} \|\nabla_K u\|_{L^2(K)} \quad \quad \text{ for } \mathcal{H}^1-\text{a.e. } x,y \in K, \label{inequS} 
\end{eqnarray}
where $dist_K(x,y)$ is the geodesic distance on $K$.  
\end{proposition}

\begin{proof} Assume first that $u \in C^\infty(\R^2)$ and let $x,y \in K$ be given. By Proposition \ref{lenghtspace} we know that there exists a geodesic Lipschitz curve $\gamma : [0,L]\to K$ with $L=\dist_K(x,y)$,  which is injective,   parametrized with constant speed so that 
$|\gamma'(t)|=1$, such that $\gamma(0)=x$, $\gamma(L)=y$ and 
$$\int_{0}^L |\gamma'(s)| \;ds= \dist_K(x,y).$$

The function $u\circ \gamma :[0,L]\to \R$ is Lipschitz continuous,  thus in particular absolutely continuous and therefore
$$u(x)-u(y)=\int_0^L \langle \nabla u\circ \gamma (t) , \gamma'(t)  \rangle \;dt,$$
from which we easily deduce that 
\begin{eqnarray}
|u(x)-u(y)| &\leq& \left(\int_0^L |\gamma'(t)|^2\, dt\right)^{\frac{1}{2}} \left(\int_0^L \langle \nabla u\circ \gamma (t) , \gamma'(t)  \rangle^2 \;dt \right)^{\frac{1}{2}} \notag \\
&\leq &\dist_K(x,y)^\frac{1}{2} \|\nabla_K u\|_{L^2(K)} , \notag
\end{eqnarray}
which  proves \eqref{inequS} in the case of a smooth $u\in \mathcal{D}(K)$. Now if $u\in H^1(K)$ we know by definition that there exists a sequence of functions $u_n \in \mathcal{D}$ such that $u_n\to u$ in $L^2(K)$ and \textcolor{black}{$\nabla_K u_n = (\nabla u_n \cdot \tau_K)\tau_K \to \nabla_K u$} in $L^2(K)$. Up to extracting a subsequence we can assume that $u_n\to u$ $\mathcal{H}^1$-a.e.\! on $K$. Applying \eqref{inequS}   to $u_n$ and then passing to the limit we  then conclude that \eqref{inequS} also holds for $u$.
\end{proof}

From Proposition \ref{embb} we get the following immediate corollary.

\begin{corollary} If $K \in \mathcal{K}(\Omega)$,  then every function $u\in H^1(K)$ admits an $L^2$-representative which is continuous.
\end{corollary}

We will also need the following $L^\infty$ estimate.

\begin{corollary}\label{bounded} If $K \in \mathcal{K}(\Omega)$,  then every function $u\in H^1(K)$ is bounded. Moreover,
$$\|u\|_{L^\infty(K)}\leq  \frac{1}{\mathcal{H}^1(K)^\frac{1}{2}}\|u\|_2+  (\mathcal{H}^1(K))^{\frac{1}{2}}\|\nabla_K u\|_2.$$
\end{corollary}

\begin{proof} For $u \in H^1(K)$ we know from Proposition \ref{embb} that for $\mathcal{H}^1$-a.e.\!  $x,y \in K$,
$$|u(x)-u(y)| \leq \|\nabla_K u\|_2\dist_K(x,y)^{\frac{1}{2}} .$$
In particular,
$$u(x)-u(y) \leq \|\nabla_K u\|_2 (\mathcal{H}^1(K))^{\frac{1}{2}},$$
thus  integrating with respect to $y \in K$  and dividing by $\mathcal{H}^1(K)$ we get 
$$u(x)-\frac{1}{\mathcal{H}^1(K)}\int_K u(y)  \;d\mathcal{H}^1(y)   \leq \|\nabla_K u\|_2 (\mathcal{H}^1(K))^{\frac{1}{2}}.$$
Finally using H\"older inequality,
$$u(x)  \leq  \|\nabla_K u\|_2 (\mathcal{H}^1(K))^{\frac{1}{2}} + \frac{1}{\mathcal{H}^1(K)^\frac{1}{2}}\|u\|_2.$$
 Reasoning the same way with $-u$ we get 
$$|u(x)|  \leq  \|\nabla_K u\|_2 (\mathcal{H}^1(K))^{\frac{1}{2}} + \frac{1}{\mathcal{H}^1(K)^\frac{1}{2}}\|u\|_2,$$
which proves the Corollary.
\end{proof}
 
We can also prove the following compact embedding result.

\begin{corollary} Let $K \in \mathcal{K}(\Omega)$.  The embedding $H^1(K)\hookrightarrow L^2(K)$ is compact. More precisely, from every bounded sequence $(u_n)_{n\in \mathbb{N}}$ in $H^1(K)$ we can extract a uniformly converging sequence, and in particular a converging sequence in $L^2(K)$.
\end{corollary}

\begin{proof} Let $(u_n)$ be a bounded sequence in $H^1(K)$. Then for each $n$ we consider the specific $L^2$ representative of $u_n$ for which  Proposition \ref{embb}  yields  the  estimate 
$$|u_n(x)-u_n(y)| \leq C\dist_K(x,y)^{\frac{1}{2}},$$
where the constant $C$ is uniform in $n$.  In particular the sequence $(u_n)$ is equicontinuous on the compact set $K\subset \R^2$.  Moreover applying Corollary \ref{bounded} we know that  $(u_n)$ is also equibounded. Thank to Arzel\`a-Ascoli theorem, we deduce the existence of a subsequence $(u_{n_k}) $ that converges uniformly on $K$. This achieves the proof of the Corollary.
\end{proof}








\section{Trace of Sobolev functions on a rectifiable set}

The object of this section is to summarize some basic facts concerning the precise representative of a Sobolev function. Consider an open set $\Omega \subset \mathbb{R}^2$. For $f\in L^1_{\rm loc}(\Omega)$, the value of the precise representative of $f$ at $x\in \Omega$ is defined by 
$$f^*(x):=\begin{cases}
	\displaystyle \lim_{r\rightarrow 0}\intbar_{B(x,r)}f(y)\,dy & \text{if the limit exists},\\
	0 & \text{otherwise}.
\end{cases}$$
The function $f^*$ depends only on the equivalence class of $f$, and coincides with $f$ a.e.\! in $\Omega$. 
In turn, we say that $f$ has an approximate limit at $x$ if there exists $t\in\mathbb{R}$ such that 
\begin{equation}\label{lebpt}
	\lim_{r\rightarrow 0}\intbar_{B(x,r)}|f(y)-t|\,dy=0\,.
\end{equation}
The set $S_f$ of points where this property fails is called the approximate discontinuity set. 

For a merely $L^1_{\rm loc}$ function, the set $S_f$ is an $\mathcal{L}^2$-negligible Borel set. On the other hand, by standard results on functions with bounded variation (see for instance~\cite[Section 3.7]{AFP2000}), we have $\mathcal{H}^{1}(S_f)=0$ whenever $f\in W^{1,1}_{\rm loc}(\Omega)$. In the sequel, any function $f$ in $W^{1,1}(\Omega)$ (or in particular in $H^1(\Omega)$) will be implicitly identified with its precise representative without mentioning it explicitly.  Since $\mathcal{H}^{1}(S_f)=0$, the pointwise values of $f$ are well defined on $K$ for $K\in \mathcal{K}(\Omega)$ and the integral 
$$\int_{K} f(x) \, d\mathcal{H}^1(x)$$
makes sense for any $f\in W^{1,1}(\Omega)$.

We  also recall the following elementary property: if $f_1\leq f_2$ a.e.\! in $\Omega$, then $f_1^*(x)\leq f_2^*(x)$ for every \textcolor{black}{$x\in \Omega\setminus (S_{f_1}\cup S_{f_2})$}. As a consequence, for any   $f\in W^{1,1}_{\rm loc}(\Omega)\cap L^\infty(\Omega)$ and  any  $K \in \mathcal{K}(\Omega)$,
\begin{equation}\label{Def:IntegraleSurKRepPrecis}
\int_{K} |f|\,d\mathcal{H}^1\leq \|f\|_{L^\infty(\Omega)} \mathcal{H}^1(K)\, .	
\end{equation}





In the sequel we would need a more accurate estimate with a Sobolev norm on the right hand side. This is the purpose of the following Lemma.

\begin{lemma}\label{Lemma:Ahlfors-bilinear-form}

Let $\Omega\subset \R^2$ be a Lipschitz domain and assume that $K\in \mathcal{K}(\Omega)$ is Ahlfors-regular, in the sense that it satisfies
\begin{equation}\label{Def:Ahlfors-constant}
\Lambda_K:=\sup\left\{\frac{\mathcal{H}^1(K\cap B(x,r))}{r}: r>0\,,\;x\in K\right\}\, <+\infty.
\end{equation}
Then for any $(v_1,v_2)\in H^1(\Omega)\times H^1(\Omega)$, the integral $\int_{K}v_1v_2\, d\mathcal{H}^1$  is well-defined. Moreover, the mapping
\[
(v_1,v_2)\in H^1(\Omega)\times H^1(\Omega)\rightarrow  \int_{K}v_1v_2\, d\mathcal{H}^1
\]
is bilinear, symmetric, nonnegative and continuous, and there exists a constant $C>0$, depending only on $\Omega$, such that
\begin{equation}\label{Estim:uv-Ahlfors}
\int_{K}|v_1v_2|\, d\mathcal{H}^1 \leq C \Lambda_K \|v_1\|_{H^1(\Omega)}\|v_2\|_{H^1(\Omega)}\, \quad \textrm{for any }(v_1,v_2)\in H^1(\Omega)\times H^1(\Omega)\, .
\end{equation}

\end{lemma}

\textcolor{black}{
\begin{remark}\label{remarkA} The standard notion of Ahlfors-regularity usually refers not only to the upper bound \eqref{Def:Ahlfors-constant}, but additionally to a lower bound 
$$\exists C_A>0 \; \text{ s.t. } \forall x\in K, \forall r\in (0,r_0), \quad  \mathcal{H}^1(K\cap B(x,r))\geq C_A r.$$
A set $K$ satisfying the above property is usually called "lower Ahlfors-regular". Here, for $K\in \mathcal{K}(\Omega)$ this property is automatically verified with $C_A=1$, since $K$ is connected. 
\end{remark}
}

\begin{remark}
	The assumption that $\Omega$ is a Lipschitz domain ensures the existence of a linear and continuous extension operator from $H^1(\Omega)$ to $H^1(\R^2)$. The same result holds under less restrictive hypotheses, for instance if we have a uniform interior cone condition (see for instance~\cite{Chenais75}).
\end{remark}

\begin{proof} The proof can be found in~\cite[Lemma 2.2]{BLM2018}, which itself relies on a standard estimate that one can find, for instance, in the book  \cite{Z1989}. For the reader's convenience, we write again the details.
	
	Define the Radon measure $\mu$ on $\R^2$ by
	\[
	\mu=\mathcal{H}^1\mres K\, .
	\]
	Notice that $\mu$ satisfies
\begin{equation}\label{Proof:Ahlfors}
	\forall x\in \R^2,\quad \forall r>0\quad \mu(B(x,r)) \leq 2r\Lambda_K\,\, .
	\end{equation}
	Indeed, let $x\in \R^2$ and $r>0$ such that $K\cap B(x,r)\neq \emptyset$. For any $z\in K\cap B(x,r)$, observe that $(K\cap B(x,r))\subset (K\cap B(z,2r))$, so that by definition~\eqref{Def:Ahlfors-constant},
	\[
	\mathcal{H}^1(K\cap B(x,r))\leq \mathcal{H}^1(K\cap B(z,2r))\leq 2r\Lambda_K\, .
	\]
	This proves~\eqref{Proof:Ahlfors}.
	
	Now, by the proof of \cite[Theorem 5.12.4]{Z1989}, there exists a universal constant $C>0$ such that
\begin{equation}\label{Resultat:Ziemer}
	 \int |w|\, d\mu  \leq C\Lambda_K \|w\|_{BV(\R^2)}\quad \textrm{for any }w\in BV(\R^2)\, ,
	\end{equation}	
	where $BV(\R^2)$ is the space of functions with bounded variation in $\R^2$. Here we have also used  that $K$ is lower Ahlfors regular with constant $1$ thanks to the connectedness assumption {\color{black}(see Remark \ref{remarkA})} in order to deduce that the constant $C$ coming from the covering Lemma \cite[Lemma 5.9.4]{Z1989}, is universal. 
	
	Let $v_1,v_2$ be in $H^1(\Omega)$. Since $\Omega$ is a Lipschitz domain, there exists a linear continuous extension operator $E:H^1(\Omega)\rightarrow H^1(\R^2)$. Now define $w:=E(v_1)E(v_2)$. The function $w$ is in $W^{1,1}(\R^2)$ so by the continuous injection $W^{1,1}(\R^2) \hookrightarrow BV(\R^2$), there exists a constant $C>0$ (depending only on $\Omega$) such that
	\begin{align*}
	\|w\|_{BV(\R^2)} & \leq C \|w\|_{W^{1,1}(\R^2)} \\
	& \leq C \|Ev_1\|_{H^1(\R^2)} \|Ev_2\|_{H^1(\R^2)} \\
	& \leq C \|v_1\|_{H^1(\Omega)} \|v_2\|_{H^1(\Omega)}\, .
	\end{align*}
Since the integral $\int |w|\, d\mu$ is equal to $\int_{K} |w|\, d\mathcal{H}^1$ and $w$ coincides with $v_1v_2$ in $\Omega$, combining the previous inequality with~\eqref{Resultat:Ziemer} yields~\eqref{Estim:uv-Ahlfors}.
\end{proof}

\section{Spectral analysis associated to the system}
\label{Spectral}
In this section we show that problem~\eqref{systeme1} has a unique global solution via variational methods. Here we follow the standard strategy of recasting the system as an abstract semilinear parabolic problem by Faedo and Galerkin (see \cite[Section 7.1.2]{Evans}).

Let $\Omega\subset \R^2$ be a bounded Lipschitz domain and $K\in \mathcal K(\Omega)$.
We start by recalling the system that we would like to solve:
\begin{equation} \label{mainSystem}
	\begin{cases}
		\partial_t v - a \Delta v = 0 &\text{in $\Omega \setminus K$}\\
		\partial_t u - b \Delta_K u = \nu v - \mu u &\text{on $K$} \\
		a \textcolor{black}{(\partial_{n}^+ + \partial_{n}^-)|_K} v = \mu u - \nu v &\text{on $K$} \\
		v = 0  &\text{on $\partial \Omega$} \\
		u = 0  &\text{on $K \cap \partial \Omega$} \\
		\partial_{n} u = 0  &\text{on $\partial K$}\setminus \partial \Omega \\
		(v,u)|_{t=0} = (v_0, u_0)
	\end{cases}
\end{equation}
 This model can be understood if $K$ is a smooth curve in $\Omega$ and $v$ and $u$ are smooth functions. In this case, the set $\partial K\setminus \partial \Omega$ is the set of endpoints of $K$ inside of $\Omega$. The Neumann condition on this set corresponds to a local mass conservation: the individuals only leave the domain through the boundary of $\Omega$.\\

We shall introduce a  weak formulation of \eqref{mainSystem} in order to make sense of the problem under fewer regularity assumptions. We take $\varphi : \Omega \to \R$ and $\psi : K \to \R$ smooth test  functions with $\varphi = 0$ on $\partial \Omega$ and $\psi = 0$ on $K\cap \partial \Omega$, we multiply the equation in $v$ by $\nu \varphi$ and the equation in $u$ by $\mu \psi$ and we integrate by parts leading us to the identities
\[
\nu \int_\Omega \partial_t v \varphi + a\nu \int_\Omega \nabla v \cdot \nabla \varphi - a \nu \int_K \partial_{n^+} v \varphi - a \nu \int_K \partial_{n^-} v \varphi = 0
\] 
and
\[
\mu \int_K \partial_t u \psi +b\mu \int_K \nabla_K u \cdot \nabla_K \psi = \mu \int_K (\nu v - \mu u) \psi.
\]
Adding the two equations together and substituting the transmission condition, we conclude that a smooth solution would solve 
\[
\nu \int_\Omega \partial_t v \varphi +  \mu \int_K \partial_t u \psi + a\nu \int_\Omega \nabla v \cdot \nabla \varphi + b\mu \int_K \nabla_K u \cdot \nabla_K \psi + \int_K (\nu v -\mu u) (\nu\varphi - \mu \psi) = 0
\] 
for any $t \in (0,T]$.

\medskip

\subsection{Variational setting} Inspired by the previous identity, we introduce some functional spaces in which to set our problem. We consider functional spaces
\[
H = H_0^1(\Omega) \times H^1_{0,\partial \Omega}(K) \quad \text{and} \quad L = L^2(\Omega) \times L^2(K)
\]
endowed with their natural scalar products and the induced topologies. More explicitly, for any $(v,u), (\varphi,\psi) \in L$, we let 
\[
\langle (v,u), (\varphi,\psi) \rangle_{L} = \nu \int_\Omega \left(v\varphi\right) dx +  \mu \int_K \left(  u \psi\right) d\mathcal{H}^1
\]
and for any $(v,u), (\varphi,\psi) \in H$,
\[
\langle (v,u), (\varphi,\psi) \rangle_H = \nu \int_\Omega \left(\nabla v \cdot \nabla \varphi + v\varphi \right) dx + \mu \int_K \left( \nabla_K  u \cdot \nabla_K \psi +  u \psi \right) d\mathcal{H}^1.
\]
It follows from Section \ref{sectionH^1} that when $K\in \mathcal{K}(\Omega)$, $H$ and $L$ are Hilbert spaces. Moreover the embedding $H \hookrightarrow L$ is compact.

We complete the variational setting by introducing the Hilbert triplet $(H^*, L = L^*, H)$, where $L^*$ and $H^*$ are the dual spaces of $L$ and $H$, respectively, and we have identified $L^*$ with $L$ itself.

We introduce the bilinear form $B : H \times H \to \R$ defined for any  $(v,u), (\varphi,\psi) \in H$ as
\[
B((v,u), (\varphi,\psi)) = a\nu \int_\Omega \nabla v \cdot \nabla \varphi \;dx + b\mu \int_K \nabla_K u \cdot \nabla_K \psi \;d\mathcal{H}^1+ \int_K (\nu v -\mu u) (\nu\varphi - \mu \psi) \;d\mathcal{H}^1.
\]

We have the following result, whose proof follows from Lemma \ref{Lemma:Ahlfors-bilinear-form}.

\begin{lemma}\label{continuousB}
	Let $K\in \mathcal{K}(\Omega)$ an  Ahlfors-regular set with upper constant $\Lambda_K$. Then the bilinear form $B$ is symmetric, continuous and coercive, that is
	\begin{itemize}
		\item for every $(v,u), (\varphi,\psi) \in H$, $B((v,u), (\varphi,\psi)) = B( (\varphi,\psi), (v,u))$;
		\item there exists $C > 0$ such that
		\[
		|B((v,u), (\varphi,\psi))| \leq C \|(v,u)\|_H \| (\varphi,\psi)\|_H \qquad \forall (v,u), (\varphi,\psi) \in H ;
		\]
		\item there exists $c > 0$ such that
		\[
		B((v,u), (v,u)) \geq c  \|(v,u)\|^2_H.
		\]
	\end{itemize} 
	Here, $c$ and $\alpha$ can be bounded below by constants that depend only on the Poincaré constant of $\Omega$ and on the Ahlfors regularity constant  $\Lambda_K$.
\end{lemma}
\begin{proof}
	The symmetry of $B$ is self evident, and the continuity follows from Lemma~\ref{Lemma:Ahlfors-bilinear-form}. Concerning the coercivity, let $C_P > 0$ and $C_T > 0$ be, respectively, the Poincaré constant of $\Omega$ and the trace constant on $K$, so that
	\[
	\int_\Omega |v|^2 dx \leq C_P \int_\Omega |\nabla v|^2 dx, \quad \int_K |v|^2 d\mathcal{H}^1 \leq C_T \int_\Omega |\nabla v|^2 dx  \qquad \forall v \in H^1_0(\Omega).
	\]
	From Lemma \ref{Lemma:Ahlfors-bilinear-form} we know that  this holds true with \textcolor{black}{$C_T=C'\Lambda_K$ where $C'>0$} is universal.

	We exploit the elementary inequality
	\begin{equation}\label{ExploitedElemIneq}
	(\alpha-\beta)^2 + \eps \alpha^2 \geq \frac{\eps}{1+\eps} \beta^2 \qquad \forall \alpha,\beta \in \R,\,  \eps \geq 0
	\end{equation}
	and find
	\[
	\begin{split}
		B((v,u), (v,u)) &=a\nu \int_\Omega |\nabla v|^2 + b\mu \int_K |\nabla_K u|^2 + \int_K |\nu v -\mu u|^2 \\
		&\geq a\nu \int_\Omega |\nabla v|^2 + b\mu \int_K |\nabla_K u|^2 -\eps \nu^2\int_K |v|^2 + \frac{\eps}{1+\eps} \mu^2 \int_K |u|^2 \\
		&\geq \frac{a\nu}{3} \int_\Omega |\nabla v|^2 + \frac{a\nu}{3C_P} \int_\Omega |v|^2 + \left( \frac{a\nu}{3C_T} - \eps \nu^2 \right)\int_K |v|^2 \\
		&\qquad + b\mu \int_K |\nabla_K u|^2 + \frac{\eps}{1+\eps} \mu^2 \int_K |u|^2.
	\end{split}
	\]
	We now choose $\eps = \frac{a}{3C_T \nu}$ and find
	\[
	B((v,u), (v,u)) \geq \min\left(\frac{a}{3}, \frac{a}{3C_P}, b, \frac{a\mu}{a + 3C_T \nu} \right) \left( \nu \int_\Omega |\nabla v|^2 + |v|^2 + \mu  \int_K |\nabla_K u|^2 + |u|^2 \right)
	\]
	for all $(v,u) \in H$. This proves the coercivity of $B$.
\end{proof}

We are now in a position to apply Fredholm's alternative to bilinear forms. Here $\delta_{hk}$ stands for Kronecker's delta symbol (that is $\delta_{hh} = 1$ and $\delta_{hk} = 0$ if $h \neq k$).
\begin{lemma}\label{Spetraum}
	Let $K\in \mathcal{K}(\Omega)$ an  Ahlfors-regular set. Then the bilinear form $B$ admits a spectral resolution: there exist sequences $\{\lambda_k\}_{k\in \N^\star} \subset [0,+\infty)$ and $\{(v_k, u_k)\}_{k\in \N^\star} \subset H$ such that
	\begin{enumerate}
		\item the sequence of \textcolor{black}{eigenvalues $\{\lambda_k\}_{k\in \N^\star}$ (counted with multiplicity) is nondecreasing} and unbounded, $\lambda_k \to +\infty$. In particular, each eigenvalue has finite multiplicity and $\lambda_1 > 0$;
		\item $\{(v_k, u_k)\}_{k\in \N^\star} \subset L$ is a orthonormal basis of $L$, and $\langle (v_k, u_k), (v_h, u_h)\rangle_L = \delta_{kh}$;
		\item $\{(v_k, u_k)\}_{k\in \N^\star} \subset H$ is a basis of $H$, and 
		\[
		B( (v_k, u_k), (v_h, u_h)) = \lambda_k \delta_{kh}.
		\]
	\end{enumerate}
\end{lemma}

\begin{proof} It is rather standard that $B$ will induce a linear operator on $L$ which has compact resolvent. Let us write some details for the convenience of the reader.  
	
	From Lemma~\ref{continuousB} we know that  $B$ is continuous and coercive on the Hilbert space $H$. By Lax-Milgram Theorem (see for instance \cite[Corollaire V.8 p.84]{brezis}) we deduce that for each $(f,g)\in L\subset H^*$, there exists a unique solution $(v,u)\in H$ to the problem 
\begin{eqnarray}
B((v,u), (\varphi,\psi))=\langle (f,g) , (\varphi,\psi) \rangle_L \quad \quad \forall \; (\varphi,\psi)\in H. \label{weakS}
\end{eqnarray}
Let $T:L\to H$ be the linear operator defined by $T(f,g)=(v,u)$, where $(v,u)$ is the solution to \eqref{weakS}. Then $T$ is a linear operator from $L\to L$ which is self-adjoint and compact. Moreover by coercivity (see Lemma \ref{continuousB}) we conclude that $\langle T(f,g) , (f,g)\rangle_L\geq c \|(v,u)\|^2_H$ for all $(f,g)\in L$.
In virtue of \cite[Theorem VI.11]{brezis} we deduce that $T$ admits an orthonormal basis of eigenfunctions \textcolor{black}{$\{(v_k, u_k)\}_{k\in \N} \subset H$} associated to a discrete sequence of eigenvalues $\mu_k > 0$,  which induces the conclusion of the Lemma by taking $\mu_k=\frac{1}{\lambda_k}$.
\end{proof}

\subsection{Lower bound on the principal eigenvalue $\lambda_1$}From the above Lemma we already know that the principal eigenvalue is positive but in the following proposition we give a more precise lower bound.

\begin{proposition} Let $K\in \mathcal{K}(\Omega)$ be  an  Ahlfors-regular set with constant $\Lambda_K$ and let us denote by $C_P$ the Poincaré constant in $\Omega$. Then there exists an explicit constant $c_0>0$ depending only on $C_P,\Lambda_K$, $a$, $\mu$, $\nu$ such that 
$$\lambda_1\geq c_0,$$
where $\lambda_1$ is the first eigenvalue  defined in Lemma \ref{Spetraum}.
\end{proposition}
\begin{proof} We denote by \textcolor{black}{$C_T=C'\Lambda_K$} the same constant as in the proof of Lemma \ref{continuousB}. First we have for $\alpha \in (0,1)$ and $\eps >0$, using once again the elementary inequality~\eqref{ExploitedElemIneq},
	\[
	\begin{split}
		B((v,u), (v,u)) &= a\nu \int_\Omega |\nabla v|^2 + b\mu \int_K |\nabla_K u|^2 + \int_K |\nu v -\mu u|^2 \\
		&\geq \frac{a\nu}{C_P} \alpha  \int_\Omega |v|^2 + \left(\frac{a\nu}{C_T} (1-\alpha) -\eps \nu^2\right)  \int_K |v|^2  + \frac{\eps}{1+\eps} \mu^2 \int_K |u|^2.
	\end{split}
	\]
	Choosing $\eps = (1-\alpha)a/(C_T\nu)$, we obtain
	\[
	B((v,u), (v,u)) \geq \frac{a\nu}{C_P} \alpha  \int_\Omega |v|^2 + \frac{(1-\alpha)a}{(1-\alpha)a+C_T \nu} \mu^2 \int_K |u|^2.
	\]	
	\textcolor{black}{Now, we seek $\alpha\in (0,1)$ such that
	\[
	\frac{a}{C_P} \alpha = \frac{(1-\alpha)a}{(1-\alpha)a+C_T \nu} \mu,
	\]
	or equivalently, $\tau(\alpha)=0$, where $\tau$ is the trinomial
	\[
	\tau(X) = a X^2 - (a+C_P \mu + C_T \nu )X + C_P\mu.
	\]
	Its discriminant is positive since $(a+C_P \mu + C_T \nu )^2-4aC_P\mu = ( a + C_T \nu - C_P \mu)^2 + 4C_TC_P\nu\mu$, hence $\tau$ possesses two real roots
	\[
	\frac{a + C_T \nu + C_P \mu}{2a} \pm \sqrt{\left( \frac{a + C_T \nu + C_P \mu}{2a} \right)^2 - \frac{C_P \mu}{a}}.
	\]
	Noticing that both roots are positive and that the product $\tau(0)\tau(1)$ is negative, the smaller root necessarily belongs to $(0,1)$. Hence, setting 
	\[
	\alpha = \frac{a + C_T \nu + C_P \mu}{2a} - \sqrt{\left( \frac{a + C_T \nu + C_P \mu}{2a} \right)^2 - \frac{C_P \mu}{a}} ,
	\] 
	we have found $\alpha\in (0,1)$ such that
}
	\begin{equation}\label{lower eigen1}
		B((v,u), (v,u)) \geq \frac{a}{C_P} \alpha \left( \nu \int_\Omega |v|^2 + \mu \int_K |u|^2\right).
	\end{equation}
	It follows that 
	\[
	\lambda_1 = \inf_{(v,u)\in H} \frac{B((v,u),(v,u))}{\|(v,u)\|_L^2} \geq \frac{a}{C_P} \alpha > 0,
	\]
	which proves the Proposition.\end{proof}

\begin{remark}
	At this stage it is not clear whether the first eigenvalue $\lambda_1$ is simple. It can be shown that the eigenfunctions associated to $\lambda_1$ do not change sign, but it is unknown if they are multiples of strictly positive functions in $\Omega \times K$ \textcolor{black}{since, at this stage, we do not know whether the equation enjoys a strong maximum principle}.
\end{remark}

\section{Existence for the Parabolic problem}
\label{Existence}
 We are now in a position to introduce the definition of a weak solution for \eqref{systeme1}.

\begin{definition}
	We consider the Hilbert space $X([0,T]; (H;H^*))$ defined as the space of functions $(v,u)$ such that
	\[
	(v,u) \in L^2([0,T]; H), \quad \frac{d}{dt} (v,u) \in L^2([0,T]; H^*)
	\]
	with its natural norm
	\[
	\|(v,u)\|_{X}^2 = \int_{[0,T]} \|(v,u)(t)\|^2_H\, dt + \int_{[0,T]} \left\|\frac{d}{dt}(v,u)(t)\right\|^2_{H^*} dt .
	\]
\end{definition}

Let us mention that by the Aubin-Lions Lemma (see \cite{LionsMagenes}), $X \subset C^0([0,T];L)$.

\begin{theorem}\label{thm main} The system \eqref{mainSystem} admits a weak solution $(v,u) \in X([0,T]; (H,H^*))$ in the following sense:  $(v,u)(0) = (v_0,u_0)$ and for every $(\varphi,\psi) \in H$ and a.e.\ $t \in (0,T]$ we have
	\begin{equation}\label{eqn weak form}
		\langle \frac{d}{dt}(v,u), (\varphi, \psi) \rangle_{(H^*,H)} + B((v,u), (\varphi,\psi)) = 0,
	\end{equation}
	where $\langle \frac{d}{dt}(v,u), (\varphi, \psi) \rangle_{(H^*,H)}$ is the duality between $H^*$ and $H$. Moreover, any weak solution is unique.
\end{theorem}

\begin{proof}  We formally write
	\[
	(v,u) = \sum_{n=1}^{+\infty} c_n(t) (v_n, u_n)
	\]
	where $c_n : [0,T] \to \R$ are continuous functions such that $c_n \in C^1((0,T];\R)$. Projecting \eqref{eqn weak form} onto the eigenfunctions we find
	\[
	\langle\sum_{m=1}^{+\infty} c_m'(t) (v_m, u_m),  (v_n, u_n) \rangle_{(H^*,H)} + B(\sum_{m=1}^{+\infty} c_m(t) (v_m, u_m),  (v_n, u_n)) = c_n'(t) + \lambda_n c_n(t) = 0
	\]
	and projecting the initial condition we get $c_n(0)= \langle (v_n, u_n),  (v_0, u_0) \rangle_{L} $.
	We therefore conclude that the couple
	\[
	(v,u) = \sum_{n=1}^{+\infty} \langle (v_n, u_n),  (v_0, u_0) \rangle_{L}  e^{-\lambda_n t} (v_n, u_n),
	\]
	is  a  weak solution, and this choice is unique.
\end{proof}

As a useful consequence we have the following.
\begin{corollary}\label{Cor:longtime}
	Under the assumptions of Theorem \ref{thm main}, assume moreover that $(v_0, u_0)$ is not orthogonal to the eigenspace of $\lambda_1$. Then there exists $(\phi, \psi) \in H$ such that
	\[
		\| (v,u) - e^{-\lambda_1 t} (\phi, \psi) \|_H = o(e^{-\lambda_1 t}) \qquad \text{as $t \to +\infty$}.	
	\]
\end{corollary}

\section{Interpretation and perspectives}\label{Sect:conclusion}

In order to interpret the previous Corollary, it seems meaningful from a modelisation perspective to compare this long-time asymptotic expansion to what would happen in the absence of a road. In that instance, one has that the solution of the problem
\begin{equation*}
	\begin{cases}
		\partial_t v - a\Delta v = 0 &\text{in $\Omega$},\\
		v = 0  &\text{on $\partial \Omega$}, \\
		v|_{t=0} = v_0&\text{in $\Omega$},
	\end{cases}
\end{equation*}
satisfies
\[
	\| v - e^{-\gamma_1 t} \phi \|_{H^1_0(\Omega)} = o(e^{-\gamma_1 t})\qquad \text{as $t \to +\infty$},	
\]
for some function $\phi \in H^1_0(\Omega)$, where $\gamma_1 > 0$ is the first eigenvalue of 
\begin{equation*}
	\begin{cases}
		 - a \Delta v_k = \gamma_k v_k &\text{in $\Omega$},\\
		v_k = 0  &\text{on $\partial \Omega$}.
	\end{cases}
\end{equation*}
As a result, to each $K\in \mathcal K(\Omega)$, one can associate the ratio $\lambda_1/\gamma_1$ which provides a criterion to determine if the presence of the road $K$ affects the effective diffusivity in the field $\Omega$, and quantify this influence.  Namely, if $\lambda_1/\gamma_1>1$, the road $K$ improves the diffusion of the population, and this effect becomes more pronounced as this ratio increases. On the opposite, for $\lambda_1/\gamma_1<1$, the road tends to slow down the spread of the species, and has no effect if $\lambda_1$ is equal to $\gamma_1$.  

\medskip

This efficiency criterion that one can associate to any admissible road $K$ belonging to the class $\mathcal K(\Omega)$ suggests several applications of the generalized framework for road-field systems introduced in this paper. The most natural one is possibly the maximization of $\lambda_1/\gamma_1$, among certain classes of roads subject to given constraints, such as:
\begin{itemize}
	\item a uniform length constraint $\mathcal H^1(K)\leq C$;
	\item topological constraints, for instance restricting the admissible class to finite unions of curves or segment;
	\item constraints involving endpoints of $K$, such as the requirement that they contain some given points in the field $\Omega$, or that some of them belonging to the boundary of $\Omega$.
\end{itemize}
Thanks to the connectedness constraint satisfied by all elements of $\mathcal K(\Omega)$, taking previous constraints into account could allow one to apply the direct method of calculus of variations and obtain compactness on any minimizing sequence in the set of compact and connected 1D-rectifiable subsets of $\Omega$, endowed with the Hausdorff distance.

\section*{Acknowledgement} A.~Zilio acknowledges support from the ANR via the project Indyana under grant agreement ANR-21-CE40-0008. The research leading to these results has received funding from the ANR project "ReaCh" (ANR-23-CE40-0023-01). M.~Bonnivard and A.~Lemenant were partially supported by the ANR project Stoiques (ANR-24-CE40-2216).


\bigskip

\bibliographystyle{alpha}
\bibliography{biblio}
\thispagestyle{empty}
\end{document}